\def\of#1{\left(#1\right)}
\def\ofc#1{\left\{#1\right\}}
\def\ofb#1{\left[#1\right]}

\def\bell{1}
\def\bellamy{2}
\def\bellamyq{3}
\def\bing{4}
\def\borsuk{5}
\def\feawri{6}
\def\fumo{7}
\def\hagopian{8}
\def\hagopianrm{9}
\def\hagopianp{10}
\def\hagopiant{11}
\def\henderson{12}
\def\iliadis{13}
\def\knaster{14}
\def\lewis{15}
\def\manka{16}
\def\scottish{17}
\def\weaklypl{18}
\def\treelike{19}
\def\perpoi{20}
\def\weakly{21}
\def\indeco{22}
\def\mitra{23}
\def\ovroa{24}
\def\ovrob{25}
\def\ovroc{26}
\def\sieklucki{27}

\input amstex

\documentstyle{amsppt}


\def\curraddr#1\endcurraddr{\address{\it Current address\/}: #1\endaddress}

\topmatter

\title A Hereditarily Indecomposable Tree-Like Continuum Without
The Fixed Point Property\\
\endtitle

\rightheadtext {A Hereditarily Indecomposable Tree-Like Continuum\dots }

\author Piotr Minc \endauthor

\address Mathematics, Auburn University, Auburn, Alabama 36849\endaddress

\email mincpio\@mail.auburn.edu \endemail
\thanks This research was supported in part by NSF grant 
\# DMS-9505108.
\endthanks

\abstract 
A hereditarily indecomposable tree-like continuum without the fixed point
property
is constructed. The example answers a question of Knaster and Bellamy.
\endabstract
\endtopmatter

\document
\baselineskip=18pt

\head
1. Introduction
\endhead
By a {\it continuum\/} we understand a nondegenerate, 
connected and compact
metric space. A continuum is {\it decomposable\/} if it
can be represented as the union of two proper subcontinua.
A continuum with no such representation is {\it indecomposable\/}.
A continuum containing no decomposable subcontinua is
called {\it hereditarily indecomposable\/}. The first example of
a hereditarily indecomposable continuum was given by B.~Knaster
in \cite{\knaster}.
A continuum is {\it tree-like\/} if it is 
the inverse limit of a sequence of trees.
This paper is motivated by the question of whether
every hereditarily
indecomposable continuum has the fixed point property.
The question was asked by Knaster (Problem 29, \cite{\lewis})
and Bellamy (page 34, \cite{\bellamyq}). In this paper we 
answer the question by
proving the following theorem. 
(Note that a fixed point is also a periodic point of period $1$.)
\proclaim{Theorem 1.1} For each positive integer $j$ there
exists a hereditarily indecomposable
tree-like continuum $X_j$
and a map $h_j:X_j\to X_j$ such that $h_j$ does not have
periodic points of periods less than or equal to $j$. 
\endproclaim

The question by Knaster and Bellamy is very closely
related to the old open problem whether 
every nonseparating
plane continuum has the fixed point property 
(see \cite{\scottish, Problem 107}).
Every nonseparating plane continuum with no interior points
is tree-like.  
In 1978, D.~P.~Bellamy \cite{\bellamy} 
presented his spectacular example of a (non planar) 
tree-like continuum 
without the fixed point property. Indecomposable continua
appear naturally in the fixed point problems for both
plane and tree-like continua.
H.~Bell (1967) \cite{\bell}, 
K.~Sieklucki (1968) \cite{\sieklucki} and 
S.~Iliadis (1970) \cite{\iliadis} proved that each fixed-point-free
map of a plane nonseparating continuum $X$ into itself must have
an invariant indecomposable subcontinuum in the boundary of $X$.
The Bellamy continuum is indecomposable.
In 1976, R.~Ma\'nka \cite{\manka} proved that every 
tree-like continuum
without the fixed point property must contain an indecomposable
continuum, not necessarily invariant (see \cite{\indeco}).
Arcs are \lq\lq{}the ultimate\rq\rq{} decomposable continua.  
In 1954, K.~Borsuk \cite{\borsuk} proved
the fixed point property for arcwise-connected tree-like continua.
The corresponding result in the plane was proven in 1971 by 
C.~Hagopian  \cite{\hagopian}. Hagopian further extended this
line of results by proving the fixed point property for
maps preserving arc-components. 
For the precise statements 
see \cite{\hagopianrm}, \cite{\hagopianp} (the plane) 
and \cite{\hagopiant} (tree-like continua).

Even though the Bellamy continuum is indecomposable,
it contains arcs playing an important role in the
construction. Very roughly speaking, the idea
is to take a map $f$ on some Knaster bucket handle continuum $K$
such that the only fixed point is the endpoint $e$ of $K$. The
map
$f$ pushes points of the arc-component of $e$ away from $e$.
Then, the fixed point $e$ is split into an infinite set $Z$
by replacing an arc containing $e$ by the cone over $Z$.
In the fixed-point-free map, branches of the cone are switched
while its vertex is pushed further away onto 
the arc-component of $e$.
Bellamy writes in \cite{\bellamyq}:
{\it It seems possible that hereditary indecomposability,
like its \lq\lq{}opposite,\rq\rq{} 
hereditary decomposability, could defeat the split,
switch and push strategy.\/} This turns out not to be the case.
In this paper we once again
(see \cite{\treelike}, \cite{\weakly} and \cite{\indeco})
start from Bellamy's marvelous construction to answer his
question. 

It was observed by J.~B.~Fugate and L.~B.~Mohler \cite{\fumo},
that if $X$ is a tree-like continuum and $f:X\to X$ is a
fixed-point-free map, then the inverse limit $\widetilde X$ of
copies $X$, with $f$ as the bonding map, is a tree-like
continuum and the shift map on $\widetilde X$ does not
have fixed points. The idea of our example is to
construct a map $f$ so that the inverse limit $\widetilde X$
is not only fixed-point-free
but it 
is also hereditarily indecomposable. 
To get our results as it is stated  in Theorem 1.1,
we start our construction from the continuum $B_j$ 
described in \cite{\treelike}.
$B_j$ is a variation of the
Bellamy example admitting a map $f_j$ with no periodic points
of periods less than or equal to $j$.
To get just a fixed-point-free map, one could begin
with the original Bellamy's continuum \cite{\bellamy}
or any of its other variations \cite{\ovroa},
\cite{\ovrob}, \cite{\ovroc} and \cite{\feawri}.
We use the theorem by Fugate and Mohler, to get
$\widetilde B_j$ and $\tilde f_j$ which is a 
\lq\lq{}smoother\rq\rq{} version of
$B_j$ and $f_j$. Every proper subcontinuum of
$\widetilde B_j$ is an arc, a
length function can be defined on arcs contained in $\widetilde B_j$ 
and $\tilde f_j$ is a homeomorphism expanding the length
of each arc.
Then, $\tilde f_j$ is corrected by a small change in such
a way that each
arc in $\widetilde B_j$ produces a 
hereditarily indecomposable continuum
in the inverse limit. For this purpose, we use very extensively
the technique presented by W.~R.~R.~Transue and the author
in \cite{\mitra} where a transitive map on $\ofb{0,1}$ 
whose inverse limit is the pseudoarc was constructed.
The first example of a map on $\ofb{0,1}$ 
whose inverse limit is the pseudoarc was given by 
G.~W.~Henderson \cite{\henderson}. 
However, Henderson's idea could not be
used in our construction, because his map is close to
the identity and we must keep close to the length expanding
map $\tilde f_j$.

The continuum constructed in this paper is, or at least
appears to be, not planar.
It would be interesting to prove that every hereditarily
indecomposable nonseparating plane continuum has the
fixed point property. Such result would parallel the
case of weakly chainable continua with a theorem
for plane nonseparating continua \cite{\weaklypl}
and a counterexample for tree-like continua 
\cite{\weakly}.

\head
2. Introducing crookedness to arcs in indecomposable continua.
\endhead
In this section, we consider a continuum $B$ satisfying the 
conditions (B1)-(B5) listed below.
We use very extensively
the technique presented by W.~R.~R.~Transue and the author
in \cite{\mitra} to get a map on $B$ whose inverse limit
is hereditarily indecomposable. The main
goal of this section is Theorem 2.13.  It corresponds very
closely to the theorem on page 1169 in \cite{\mitra}
and we have to virtually repeat the same argument to see
whether its fine details work in the new setting.

Suppose that $B$ is a tree-like continuum and 
$\beta$ is a map of $B$ onto $\ofb{0,1}$ such that the following
conditions are satisfied:
\roster
\item"{(B1)}" each (non-trivial) proper subcontinuum of 
$B$ is an arc,
\item"{(B2)}" $\dim\of{\beta^{-1}\of t}=0$ for each $t\in\ofb{0,1}$, 
\item"{(B3)}" $\beta^{-1}\of{\of{0,1}}$ 
is homeomorphic to the product
$\of{0,1}$ and the Cantor set, 
\item"{(B4)}" $\beta$ restricted to $\beta^{-1}\of{\of{0,1}}$ is
the projection of the product onto the $\of{0,1}$ component, and
\item"{(B5)}" every point of $B$ is either an endpoint or it has
a neighborhood homeomorphic to the product of the Cantor set and
an open interval. (A point of $B$ is an {\it endpoint\/} if it does
not belong to an open arc contained in $B$.) 
\endroster
For an arc
$A\subset B$, we will define the {\it $\beta$-length\/} of $A$
as the sum
$$\lambda\of A=\sum_{C\in\Cal C\of A} 
\operatorname{diam}\of{\beta\of C}$$
where $\Cal C\of A$ denote the set of
components of
$A\cap \beta^{-1}\of{\of{0,1}}$.
Observe that $\lambda\of{A}=\lambda\of{A_1}+\lambda\of{A_2}$
if $A_1$ and $A_2$ are subarcs of $A$ such that
$A=A_1\cup A_2$ and $A=A_1\cap A_2$ is a single point.

We will leave the proof of 
the following proposition to the reader.
\proclaim{Proposition 2.1} For each positive number $\eta$
there is a positive number $\epsilon$ with the property
that $\operatorname{diam}\of A<\eta$ for each arc 
$A\subset B$ such that $\lambda\of A<\epsilon$.
\endproclaim

For a point $a\in\beta^{-1}\of{\of{0,1}}$, let $L\of a$ denote
the closure of the component of $a$ in $\beta^{-1}\of{\of{0,1}}$.
Observe that $\beta$ restricted to $L\of a$ is a homeomorphism
onto $\ofb{0,1}$. Let $h\ofb a:\ofb{0,1}\to L\of a$ be the homeomorphism
so that $\beta\circ h\ofb a$ is the identity on 
$\ofb{0,1}$. Let $e\of a$ and $d\of a$ denote the endpoints of
$L\of a$ so that $\beta\of{e\of a}=0$ and $\beta\of{d\of a}=1$.

Let $Z=\beta^{-1}\of{\frac12}$. 
For each point $z\in Z$ there is
at most one point $\tau^-\of z\in Z\setminus\ofc z$ such that
$e\of{\tau^-\of z}=e\of z$. In the case when 
$e\of z \ne e\of t$ for each point
$t\in Z\setminus\ofc z$, we will set $\tau^-\of z =z$. Similarly,
there is
at most one point $\tau^+\of z\in Z\setminus\ofc z$ such that
$d\of{\tau^+\of z}=d\of z$. In the case when 
$d\of z \ne d\of t$ for each point
$t\in Z\setminus\ofc z$, we will set $\tau^+\of z =z$. 

\proclaim{Proposition 2.2}
$\lim_{i\to\infty}\tau^-\of{z_i}=\tau^-\of{z}$ and
$\lim_{i\to\infty}\tau^+\of{z_i}=\tau^+\of{z}$ for each
sequence $z_i\in Z$ such that $\lim_{i\to\infty}z_i=z$.
\endproclaim
\demo{Proof}Suppose $s$ is the limit of a subsequence
of $\of{\tau^-\of{z_i}}_{i=0}^\infty$. Observe that 
$e\of s=e\of z$. Thus, either $s=z$ or $s=\tau^-\of z$.
Suppose $s=z$. In this case, an infinite sequence
sets of the form $L\of{z_i}\cup L\of{\tau^-\of{z_i}}$ 
converges $L\of z$. If $\tau^-\of{z_i}=z_i$, $e\of{z_i}$
is an endpoint. If $\tau^-\of{z_i}\ne z_i$, the set
$L\of{z_i}\cup L\of{\tau^-\of{z_i}}$ forms a "V" close
to $L\of z$. Since there is either infinitely many cases
of $\tau^-\of{z_i}=z_i$ or infinitely many  cases of
$\tau^-\of{z_i}\ne z_i$, the point $e\of z$ does not have
a neighborhood homeomorphic to the product of the Cantor set and
an open interval. So, by (B5), $e\of z$ is an endpoint. Thus,
$\tau^-\of z=z$ and consequently $s=\tau^-\of z$.
Since the limit of an arbitrary subsequence
of $\of{\tau^-\of{z_i}}_{i=0}^\infty$ is $\tau^-\of z$, we have the
result that $\lim_{i\to\infty}\tau^-\of{z_i}=\tau^-\of{z}$.
The proof of
$\lim_{i\to\infty}\tau^+\of{z_i}=\tau^+\of{z}$ is similar 
and will be omitted.
\enddemo

Let $f:B\to B$ be a map.
We will say that $f$ is {\it Lipschitz\/} if there is a constant
$s$ such that $\lambda\of{f\of A}<s\lambda\of A$ for each arc 
$A\subset B$. If $\sigma$ is a number greater than $1$, we will 
say that 
$f$ {\it stretches\/} by $\sigma$ if
$\lambda\of{f\of A}\ge\sigma\lambda\of A$ for each arc 
$A\subset B$.

For any two points $a$ and $b$ belonging 
to the same arc-component of
$B$, let $\ofb{a,b}$ denote the arc in $B$ between $a$ and $b$.
($\ofb{a,a}$ is just a single point.)

Suppose $f_1$ and $f_2$ are maps of $B$ into itself such that
$f_1\of t$ and $f_2\of t$ belong to the same arc component of
$B$ for each $t\in B$. Let $d_\lambda\of{f_1,f_2}$ be the
supremum of $\lambda\of{\ofb{f_1\of t,f_2\of t}}$ where
$t\in B$. By saying that $d_\lambda\of{f,f^\prime}$ is finite
we will imply, in particular, that $f\of t$ and $f^\prime\of t$ 
belong to the same arc component of
$B$ for each $t\in B$.

The following proposition is a simple consequence of
Proposition 2.1. 
\proclaim{Proposition 2.3} For each positive number $\eta$
there is a positive number $\epsilon$ with the property
that the distance between $f$ and $f^\prime$ is less than 
$\eta$ for every two maps $f$ and $f^\prime$ of $B$ 
into itself such that
$d_\lambda\of{f,f^\prime}<\epsilon$.
\endproclaim

\proclaim{Proposition 2.4}Suppose $f:B\to B$ is a Lipschitz
map with the Lipschitz constant $s$. Let $F:B\to B$ be a map
such that $d_\lambda\of{f,F}<\eta$ for a certain positive 
number $\eta$. Then $d_\lambda\of{f^j,F^j}<\of{s+1}^{j-1}\eta$
for each positive integer $j$.
\endproclaim
\demo{Proof} The proposition is true for $j=1$. We assume that
it is true for $j=k$ and we will prove it for $j=k+1$.
Let $t$ be an arbitrary point of $B$. Observe that
$$\lambda\of{\ofb{f^{k+1}\of t,f^k\of{F\of t}}}<
s^k\lambda\of{\ofb{f\of t,F\of t}}<s^k\eta.$$ Since
$$\lambda\of{\ofb{f^{k+1}\of t,F^{k+1}\of t}}\le
\lambda\of{\ofb{f^{k+1}\of t,f^k\of{F\of t}}}+
\lambda\of{\ofb{f^k\of{F\of t},F^k\of{F\of t}}},$$
we have the result that
$$\lambda\of{\ofb{f^{k+1}\of t,F^{k+1}\of t}}<
s^k\eta+\of{s+1}^{k-1}\eta<\of{s+1}^k\eta
$$ and the proposition is proven by induction.
\enddemo

For any non-degenerate arc $A\subset B$, one of its endpoints $a$
and any positive number $\epsilon$, 
let $H\of{a,A,\epsilon}$ denote the set of points $t\in B$ such that
$a$ and $t$ belong to the same arc component of $B$, 
$\ofb{a,t}\cap A=\ofc a$ and
$\lambda\of{\ofb{a,t}}\le\epsilon$. 
Clearly, 
$H\of{a,A,\epsilon}$ is either the point $a$ or
an arc with $a$ as one of its endpoints.

For any arc $A\subset B$, suppose that $a$ and $b$ denote the 
endpoints of 
$A$. If $\epsilon$ is a positive number, by
$N\of{A,\epsilon}$ we will denote the union 
$H\of{a,A,\epsilon}\cup A\cup H\of{b,A,\epsilon}$.

Observe that the following proposition holds.
\proclaim{Proposition 2.5} Suppose $f:B\to B$ is a Lipschitz
map with the Lipschitz constant $s$. Then
$f\of{N\of{A,\epsilon}}\subset N\of{f\of A, s\epsilon}$
for any arc $A\subset B$ and any positive number
$\epsilon$.
\endproclaim

For any arc $A\subset B$ and a positive number $\epsilon$ such that
$\epsilon<\lambda\of A/2$, let
$K_1\of{A,\epsilon}$ be a subarc of $A$ containing one of
the endpoints of $A$
so that $\lambda\of{K_1\of{A,\epsilon}}=\lambda\of A-\epsilon$.
Let 
$K_2\of{A,\epsilon}$ be the subarc of $A$ containing 
the endpoint of $A$ that does not belong to $K_1\of{A,\epsilon}$
so that $\lambda\of{K_2\of{A,\epsilon}}=\lambda\of A-\epsilon$.
Additionally, let
$K\of{A,\epsilon}=K_1\of{A,\epsilon}\cap K_2\of{A,\epsilon}$.

Let $f:B\to B$ be a map and let $\epsilon$ and $\delta$ be two
positive numbers. We will say that $f$ is 
$\of{\epsilon,\delta}$-{\it crooked\/} if for each arcs
$A$ and $C$ contained in $B$ so that
$2\delta<\lambda\of A<\epsilon$ and
$f\of C=A$, there are two disjoint arcs $C_1$ and $C_2$ contained
in $C$ so that $f\of{C_1}=f\of{C_2}=K\of{A,\delta}$.

The proof of the next three propositions is left to the reader.
\proclaim{Proposition 2.6}Let $\epsilon$ and $\delta$ be two
positive numbers and let $f:B\to B$ be a 
$\of{\epsilon,\delta}$-crooked map.
Then, 
for each arcs
$A$ and $C$ contained in $B$ so that
$2\delta<\lambda\of A<\epsilon$ and
$f\of C=A$, there are two disjoint arcs $C_1^\prime$ and 
$C_2^\prime$ contained
in $C$ so that $f\of{C_1^\prime}=K_1\of{A,\delta}$
and $f\of{C_1^\prime}=K_2\of{A,\delta}$.
\endproclaim
\proclaim{Proposition 2.7}If
$\alpha\ge\delta$ and $\mu\le\epsilon$, 
then every $\of{\epsilon,\delta}$-crooked
map of $B$ into itself is also $\of{\mu,\alpha}$-crooked.
\endproclaim

\proclaim{Proposition 2.8} Let $\epsilon$, $\delta$ and
$\eta$ be positive numbers.
Suppose $f$ and $F$ are maps
of $B$ into itself such that
$d_\lambda\of{f,F}<\eta$.
If $f$ is $\of{\epsilon,\delta}$-crooked, then
$F$ is $\of{\epsilon,\delta+2\eta}$-crooked.
\endproclaim

The following proposition corresponds 
to \cite{\mitra, Proposition 4}.
\proclaim{Proposition 2.9}Let $f$ be a
map
of $B$ into itself with the property that
for each positive numbers $\mu$ and $\delta$,
there is a positive integer $n$ such that
$f^n$ is
$\of{\mu,\delta}$-crooked.
Then the inverse limit of copies of $B$ with $f$
as the bonding map is hereditarily indecomposable.
\endproclaim
\demo{Proof}
Let $X$ denote the inverse limit of copies of $B$ with $f$
as the bonding map. Let $p_i$ denote the projection into 
$i$-th copy of $B$ in the inverse sequence. Suppose that
$X$ contains two subcontinua $X_1$ and $X_2$ such that
the sets $X_1\setminus X_2$, $X_2\setminus X_1$ and
$X_1\cap X_2$ are not empty. There is a positive integer $k$
such that the sets $p_k\of{X_1}\setminus p_k\of{X_2}$ and
$p_k\of{X_2}\setminus p_k\of{X_1}$ are not empty.
By (B1), $p_k\of{X_1}$ and $p_k\of{X_2}$ are arcs.
It follows that $A=p_k\of{X_1}\cup p_k\of{X_2}$ is also an
arc. Take $\mu>\lambda\of{A}$.
Let $A_1$  denote the closure of 
 $p_k\of{X_1}\setminus p_k\of{X_2}$ and
let 
$A_2$ be the closure of 
 $p_k\of{X_2}\setminus p_k\of{X_1}$. Clearly, $A_1$ and $A_2$
are nondegenerate arcs. Let $\delta$ be a positive number
less than the minimum of $\lambda\of{A_1}$ and $\lambda\of{A_2}$.

There is a positive integer $n$ such that
$f^n$ is
$\of{\mu,\delta}$-crooked.
Let $C=p_{k+n}\of{X_1}\cup p_{k+n}\of{X_2}$. Since
$f^n\of C=A$, $C$ is a proper subcontinuum of $B$. Consequently,
$C$ is an arc. There are two disjoint arcs $C_1$ and $C_2$ 
contained in $C$ such that $f^n\of{C_1}=f^n\of{C_2}=K\of{A,\delta}$.
Since $K\of{A,\delta}$ is contained in neither 
$p_k\of{X_1}$ nor
$p_k\of{X_2}$, we have the result that, for $i=1,2$, 
$C_i$ is  contained in neither 
$p_{k+n}\of{X_1}$ nor
$p_{k+n}\of{X_2}$. It follows that each of the arcs
$C_1$ and $C_2$ contains $p_{k+n}\of{X_1}\cap p_{k+n}\of{X_2}$,
and thus $C_1\cap C_2\ne\emptyset$, a contradiction. 
\enddemo

The next proposition readily follows from the definitions
and Proposition 2.3.
\proclaim{Proposition 2.10}Let
$f_0,f_1,f_2,\dots$ be a sequence of maps of
$B$ into itself
such that $\sum_{i=0}^\infty d_\lambda\of{f_i,f_{i+1}}$ is
finite. Then the sequence $f_0,f_1,f_2,\dots$ 
converges uniformly and its limit $f$ has the property that
$d_\lambda\of{f_0,f}\le 
\sum_{i=0}^\infty d_\lambda\of{f_i,f_{i+1}}$.
Additionally, if each of the maps $f_0,f_1,f_2,\dots$ is
$\of{\mu,\delta}$-crooked for some positive numbers
$\mu$ and $\delta$, then $f$ is also $\of{\mu,\delta}$-crooked.
\endproclaim

The following proposition corresponds 
to \cite{\mitra, Proposition 5}.
\proclaim{Proposition 2.11}
Let $\epsilon<1$ and $\gamma<\epsilon/4$ be two positive numbers.
Then there is a Lipschitz map $g:B\to B$ such that
\roster
\item"{(i)}" for each point $t\in B$, $t$ and $g\of t$ 
are in the same
arc component of $B$ and
$\lambda\of{\ofb{t,g\of t}}<\epsilon/2+\gamma$,
\item"{(ii)}" $g$ is $\of{\epsilon,\gamma}$-crooked, 
\item"{(iii)}" $\lambda\of{g\of A}\ge \lambda\of A$ for each
arc $A\subset B$,
\endroster
and if, additionally, $\lambda\of A\ge\gamma$, then
\roster
\item"{(iv)}" $\lambda\of{g\of A}>\epsilon/2$, and
\item"{(v)}" $g\of{N\of{A,r}}\subset N\of{g\of A,r+\gamma}$ for
each positive number  $r$. 
\endroster
\endproclaim
\demo{Proof}
Let $q$ be an integer such that $q\gamma\ge4$.
Let $\mu=1/q$ and let $p$ be the integer so that
$\mu\of{p-2}<\epsilon/2\le\mu\of{p-1}$.
We will now define a piece-wise linear function
$g_0:\ofb{0,1}\to\ofb{-1,2}$. 
For each $i=0,\dots,q-1$, on the interval of
the form
$\ofb{i\mu,\of{i+1}\mu}$ the function $g_0$ is defined
in the following way: 
$g_0\of{i\mu}=i\mu$, then we move up to $\of{i+p}\mu$ 
in crooked fashion, joining linearly values 
$\,i\mu$, $\,\of{i+1}\mu$, $\,\of{i+2}\mu$, $\,\of{i+1}\mu$,
$\,\of{i+2}\mu$, $\,\of{i+3}\mu$, $\,\of{i+2}\mu$,
$\,\of{i+1}\mu$, $\,\of{i+2}\mu$, $\,\of{i+3}\mu$, 
$\,\of{i+2}\mu$, $\,\of{i+3}\mu$, $\,\of{i+4}\mu$, ..., 
(see\cite{\bing}), then we move crookedly down to 
$\of{i-p}\mu$, and then up again to $\of{i+1}\mu$,
a goal which is achieved at $\of{i+1}\mu$.

Now, we are ready to define $g$. Let $z$ be an arbitrary point
of $Z$ and let $t$ be an arbitrary point of $L\of z$.
Let
$$g\of t=\cases
h\ofb{z}\of{g_0\of t},&\text{if $0\le g_0\of t\le1$}\\
h\ofb{\tau^+\of z}\of{2-g_0\of t},&\text{if $g_0\of t>1$}\\
h\ofb{\tau^-\of z}\of{-g_0\of t},&\text{if $0> g_0\of t$}.
\endcases
$$ 
Continuity of $g$ is guarantied by Proposition 2.2. 
It may be verified
that $g$ has the remaining required properties.
\enddemo

The following Lemma 2.12 corresponds 
to the lemma on page 1167 in \cite{\mitra}.
The proof presented here is almost the same as in \cite{\mitra}.
The only difference, apart from setting it for the
continuum $B$ instead of the interval $\ofb{0,1}$, is replacing
the condition that $f$ eventually expands each interval
to $\ofb{0,1}$ by the condition that $f$ stretches.
\proclaim{Lemma 2.12} Let $\sigma$ be a number greater than $1$
and let $f:B\to B$ be a Lipschitz map that stretches by $\sigma$.
Let $\eta$, $\delta$ and $\mu$ be positive numbers. Then there is
a map $F:B\to B$ and there is a positive integer $n$ such that
\roster
\item"{}"
$F$ is Lipschitz and stretches by $\sigma$, 
\item"{}" $d_{\lambda}\of{f,F}<\eta$, and
\item"{}" $F^n$ is $\of{\mu,\delta}$-crooked.
\endroster
\endproclaim
\demo{Proof}
Since $f$ is Lipschitz, there is a number $s>2$ such that
$$\lambda\of{f\of A}<s\lambda\of A \text{for each arc $A\subset B$.}
\tag1$$
Set $\epsilon=\eta/s$. Let $n$ be a positive integer such that
$\sigma^n\epsilon>2\mu$. Let $\gamma$ be a positive number less than
$\min\of{\epsilon/4,\delta s^{-n}/4}$.

Let $g$ be a map satisfying Proposition 2.11. Define
$F=f\circ g$. Clearly, $F$ is Lipschitz as the composition of two
Lipschitz maps. Since $f$ stretches by $\sigma$, it follows from
Proposition 2.11 (iii) that
$F$ also stretches by $\sigma$.

Suppose $t$ is an arbitrary point of $B$. By Proposition 2.11 (i),
$t$ and $g\of t$ are in the same arc component of $B$ and
$\lambda\of{\ofb{g\of t,t}}<\epsilon=\eta/s$. It follows
that $F\of t$ and $f\of t$
are in the same arc component of $B$ and,
by \thetag{1},
$\lambda\of{\ofb{F\of t, f\of t}}<\eta$. 

To prove that $F^n$ is $\of{\mu,\delta}$-crooked, we 
need the following two claims.

\proclaim{Claim 1} Let $A\subset B$ be an arc such that
$\lambda\of A\ge\gamma$. Let $r$ be a positive real number
and let $j$ be a positive integer.
Then
$F^j\of{N\of{A,r}}\subset N\of{F^j\of A,q}$,
where $q=s^j\of{r+2\gamma}$.
\endproclaim
\demo{Proof of Claim {\rm 1}}
Since $F$ stretches by $\sigma>1$, we have the result that
$\lambda\of{F^i\of A}>\gamma$ for each positive integer
$i$. Set $q_0=r$ and $q_{i+1}=s\of{q_i+\gamma}$ 
for $i=0,1,\dots$. Repeatedly using Proposition 2.11 (v),
the condition
 \thetag{1} and Proposition 2.5, one can prove
$F^j\of{N\of{A,r}}\subset N\of{F^j\of A,q_j}$.
Observe that
$q_j=s^jr+s^j\gamma+s^{j-1}\gamma+\dots+s\gamma=
s^jr+\gamma s\frac{s^j-1}{s-1}=
s^j\of{r+\gamma\frac s{s-1}\frac{s^j-1}{s^j}}<
s^j\of{r+\gamma\frac s{s-1}}$.
Since $s>2$, $\frac s{s-1}<2$, and, consequently,
$q_j<s^j\of{r+2\gamma}$, so the claim is true.
\enddemo

\proclaim{Claim 2} Let $C\subset B$ be an arc such that
$2\gamma<\lambda\of C<\epsilon +2\nu$ for some positive
number $\nu$. Then there is an arc $G\subset C$
such that $2\gamma<\lambda\of G<\epsilon$ and
$C\subset N\of{G,\nu}$.
\endproclaim
\demo{Proof of Claim {\rm 2}}
If 
$\lambda\of{C}<\epsilon$, we set $G=C$. So,
we may assume that $\lambda\of{C}\ge\epsilon$.
Since
$\lambda\of{C}-\epsilon<2\nu$,
there is a number $\kappa>\lambda\of{C}-\epsilon$
such that
$\kappa<\lambda\of{C}-\epsilon +\gamma$
and
$\kappa<2\nu$.
Set $G=K\of{C,\kappa/2}$. In this case, 
$\lambda\of G=\lambda\of{C}-\kappa<\epsilon$. We have also
$\lambda\of G=\lambda\of{C}-\kappa>\epsilon-\gamma>2\gamma$.
Finally,
$C=N\of{G,\kappa/2}\subset N\of{G,\nu}$.
\enddemo

Let $A$ be a subarc of $B$ such that 
$2\delta<\lambda\of A<\mu$. Suppose that
$I\subset B$ is an arc such that $F^n\of I=A$.
To complete the proof of the lemma, it is
enough to show that
$$\hbox{\text{\vbox{\hsize 3.5in\noindent
there are two disjoint subarcs
$A_1$ and $A_2$ of $I$ such that
$K\of{A,\delta}\subset F^n\of{A_1}
\cap F^n\of{A_2}$.}}}
\tag2
$$

We will observe that 
$$\lambda\of I<\gamma.\tag3$$
Suppose to the contrary. Then, by Proposition 2.11 (iv),
$\lambda\of{g\of I}>\epsilon/2$. Since $f$ stretches
by $\sigma$, $\lambda\of{F\of I}>\sigma\epsilon/2$.
Since $F$ stretches by $\sigma$,
$\lambda\of{F^n\of I}>\sigma^n\epsilon/2$.
The last number is greater than $\mu$ by the choice of $n$,
so we get $\lambda\of A>\mu$, a contradiction.

Since $F^n\of I=A$ and $\lambda\of A>2\delta>\gamma$,
$\lambda\of{F^n\of I}>\gamma$. Let $m$ be the greatest
integer such that $\lambda\of{F^m\of I}<\gamma$. Clearly,
$0\le m<n$. Denote $F^m\of I$ by $M$.
We will consider two cases: 
$\lambda\of{g\of M}>2\gamma$ and 
$\lambda\of{g\of M}\le2\gamma$.

\proclaim\nofrills{Case:}\quad$\lambda\of{g\of M}>2\gamma$.
\endproclaim
\noindent 
Since $\lambda\of M<\gamma$, it follows from Proposition 2.11 (i)
that $\lambda\of{g\of M}<\epsilon+3\gamma$.
It follows from Claim 2 that there is an arc $G\subset g\of M$
such that $2\gamma<\lambda\of G<\epsilon$ and
$g\of M\subset N\of{G,2\gamma}$. 
Since $g$ is $\of{\epsilon,\gamma}$-crooked 
(Proposition 2.11 (ii)),
it follows from Proposition 2.6 that there two disjoint
arcs $M_1$ and $M_2$ contained in $M$ such that
$g\of{M_1}=K_1\of{G,\gamma}$ and
$g\of{M_2}=K_2\of{G,\gamma}$.
There are two disjoint arcs $A_1$ and $A_2$ contained
in $I$ such that $F^m\of{A_k}=M_k$ for $k=1,2$.
Observe that $g\of M\subset N\of{g\of{M_k},3\gamma}$
and it follows from Proposition 2.5 that
$F\of M\subset N\of{F\of{M_k},3s\gamma}$ for $k=1,2$.
Since $\lambda\of{g\of{M_k}}>\gamma$ and $f$ stretches
(by $\sigma$), we have the result that
$\lambda\of{F\of{M_k}}>\gamma$. It follows from Claim 1
that
$A=F^n\of I=F^{n-m}\of M\subset
N\of{F^{n-m}\of{M_k},s^{n-m-1}\of{3s\gamma+2\gamma}}$.
Since 
$s^{n-m-1}\of{3s\gamma+2\gamma}<4s^n\gamma<\delta$,
$K\of{A,\delta}\subset 
F^{n-m}\of{M_1}\cap F^{n-m}\of{M_2}=
F^{n}\of{A_1}\cap F^{n}\of{A_2}$.

\proclaim\nofrills{Case:}\quad$\lambda\of{g\of M}\le2\gamma$.
\endproclaim
\noindent It follows from \thetag{1} and the choice of $m$ that
$\gamma\le\lambda\of{F\of M}<2\gamma s$. (Notice, that
the last inequality implies, in particular, that
$n>m+1$, because $\gamma s<\delta$.)

By Proposition 2.11 (i) and (iv), 
$\epsilon/2<\lambda\of{g\of{F\of M}}<\epsilon+2\gamma\of{s+1}$.
Using Claim 2 we get an arc $G\subset g\of{F\of M}$
such that $2\gamma<\lambda\of G<\epsilon$ and
$g\of{F\of M}\subset N\of{G,\gamma\of{s+1}}$.
Since $g$ is $\of{\epsilon,\gamma}$-crooked (Proposition 2.11 (ii)),
it follows from Proposition 2.6 that there two disjoint
arcs $M_1$ and $M_2$ contained in $F\of M$ such that
$g\of{M_1}=K_1\of{G,\gamma}$ and
$g\of{M_2}=K_2\of{G,\gamma}$.
There are two disjoint arcs $A_1$ and $A_2$ contained
in $I$ such that $F^{m+1}\of{A_k}=M_k$ for $k=1,2$.
Observe that $g\of{F\of M}\subset 
N\of{g\of{M_k},\gamma\of{s+1}+\gamma}\subset
N\of{g\of{M_k},2s\gamma}
$
and it follows from Proposition 2.5 that
$F^2\of M\subset N\of{F\of{M_k},2s^2\gamma}$ for $k=1,2$.
Since $\lambda\of{g\of{M_k}}>\gamma$ and $f$ stretches
(by $\sigma$), we have the result that
$\lambda\of{F\of{M_k}}>\gamma$. 
It follows from Claim 1
that
$A=F^n\of I=F^{n-m}\of M\subset
N\of{F^{n-m-1}\of{M_k},s^{n-m-2}\of{2s^2\gamma+2\gamma}}$.
Since 
$s^{n-m-2}\of{2s^2\gamma+2\gamma}
<4s^n\gamma<\delta$,
$K\of{A,\delta}\subset 
F^{n-m}\of{M_1}\cap F^{n-m}\of{M_2}=
F^{n}\of{A_1}\cap F^{n}\of{A_2}$ and the proof of the lemma
is complete.
\enddemo
\proclaim{Theorem 2.13} Suppose $\epsilon$ is a positive number
and $f_0:B\to B$ is a Lipschitz map
that stretches by a certain number $\sigma>1$.
Then there is a map $f:B\to B$ such that
$d_\lambda\of{f,f_0}<\epsilon$ and the inverse limit
of copies of $B$ with $f$ as the bonding map is 
hereditarily indecomposable.
\endproclaim
\demo{Proof}
We are going to construct a sequence of positive integers
$n\of1$, $n\of2$, $n\of3$, $\dots$ and a sequence $f_1,f_2,f_3,\dots$ of
maps of $B$ into itself such that,
for each positive integer $i$,
the following conditions are satisfied:
\roster
\item"{(i)}" $f_i$ is a Lipschitz map stretching by $\sigma$,
\item"{(ii)}" $d_\lambda\of{f_{i-1},f_i}<2^{-i}\epsilon$, and
\item"{(iii)}" $f_i^{n\of k}$ is $\of{k,2^{-k}-2^{-k-i}}$-crooked
for each positive integer $k\le i$.
\endroster
To construct $f_1$, we apply Lemma 2.12 with
$f=f_0$, $\eta=\epsilon/2$, $\delta=2^{-2}$ and $\mu=1$.
We set $f_1=F$ and $n\of1=n$. We assume that $n\of1,n\of2,\dots,
n\of{i-1}$ and $f_1,f_2,\dots,f_{i-1}$ have already been constructed,
and we will construct $n\of i$ and $f_i$.

By Proposition 2.4, there is a positive number $\eta<2^{-i}\epsilon$
such that, if $F:B\to B$ is a map with the property 
$d_\lambda\of{f_{i-1},F}<\eta$, then
$d_\lambda\of{f_{i-1}^{n\of k},F^{n\of k}}<2^{-k-i-1}$ for
each positive integer $k<i$. 
Now, use Lemma 2.12 with $f=f_{i-1}$, $\delta=2^{-i}-2^{-i-i}$
and $\mu=i$. Define $f_i=F$ and $n\of i=n$.
Clearly, the conditions (i), (ii) and (iii) for $k=i$ are satisfied.
By Proposition 2.8, the choice of
$\eta$ guarantees that $f_i^{n\of k}=F^{n\of k}$ is
$\of{k,2^{-k}-2^{-k-i}}$-crooked for $k<i$. So the construction is
complete.

By (ii) and Proposition 2.10, the sequence $f_0,f_1,f_2,\dots$ 
converges uniformly. Denote the limit by $f$.
By Proposition 2.10, $d_\lambda\of{f,f_0}<\epsilon$.
By (iii) and Proposition 2.7, 
$f_i^{n\of k}$ is $\of{k,2^{-k}}$-crooked for every
positive integers $k$ and $i$ such that $k\le i$.
By the second part of proposition 2.10, 
$f^{n\of k}$ is $\of{k,2^{-k}}$-crooked for every
positive integer $k$. Applying propositions 2.7 and 2.9, we
get the result that
the inverse limit
of copies of $B$ with $f$ as the bonding map is 
hereditarily indecomposable.
\enddemo

\head
3. The main result
\endhead
In this section we finish our proof of Theorem 1.1
by supplying a suitable continuum $B$ to the machinery
developed in the previous section.
For this purpose we will take $\widetilde B_j$ as
it is described in \cite{\perpoi}.
The continuum $\widetilde B_j$ results from applying
the theorem by Fugate and Mohler \cite{\fumo} to
the continuum $B_j$ described in \cite{\treelike}.
As we noted in the introduction,
to get just a fixed-point-free map in Theorem 1.1, 
one could replace $B_j$
with the original Bellamy's continuum \cite{\bellamy}
or any of its other variations \cite{\ovroa},
\cite{\ovrob}, \cite{\ovroc} and \cite{\feawri}.
Even though the properties required in
Section 2 are quite apparent for
any of
these continua, we feel obliged to check the details
at least in the case of $B_j$. To do that we need to
summarize the construction from \cite{\treelike} and
Section 2 of \cite{\perpoi}.

For each positive integer $k$, let $g_k:\ofb{0,1}\to \ofb{0,1}$ 
be the function stretching the interval $\ofb{0,1}$ $k$ times
and then folding it uniformly back onto itself. For example,
$g_2$ is the roof-top map on $\ofb{0,1}$.
For each positive integer $n$,  let $S_n$ be the inverse limit of the 
inverse system of copies of $\ofb{0,1}$ with every bonding map equal to 
$g_n$.
Let $p_n^k$ be the projection of $S_n$ onto the $k$-th element of the
inverse system. Let $e_n$ denote the point 
$\left(0,0,\dots\right)$ and let
$d_n$ denote the point $\left(1,1/n,1/n^2,1/n^3,\dots\right)$. 
Let $J_n$ denote the arc
in $S_n$ between $e_n$ and $d_n$. Let $g$ denote the map from $S_n$
onto itself induced by $g_2$, i.e. 
$g\of{\of{x_0,x_1,\dots}}=\of{
g_2\of{x_0},g_2\of{x_1},\dots}$.

The following properties of $S_n$ are very well known.
\roster
\item"{(S1)}" Every proper subcontinuum of $S_n$ is an arc.
\item"{(S2)}" For each arc $A\subset S_n$, there is a positive
integer $m_0$ such that $g^m\of{d_n}\notin A$ for each
$m\ge m_0$.
\endroster
If $n$ is even, we additionally have:
\roster
\item"{(S3)}" $e_n$ is the only endpoint of $S_n$,
\item"{(S4)}" each point of $S_n\setminus\ofc{e_n}$
has a neighborhood homeomorphic to the product of
the Cantor set and an open interval,
\item"{(S5)}" $g$ is a homeomorphism, and
\item"{(S6)}" $g^{-1}\of{J_n}\subset J_n\setminus\ofc{d_n}$.
\endroster

For each positive integer $j$,  
$n\of j$ denotes
$2\of{4^1-1}\of{4^2-1}\dots\of{4^j-1}$. 
By a slight variation of the original Bellamy's construction 
\cite{\bellamy}, 
it was
proven in \cite{\treelike}
 that there is a tree-like continuum $B_j$ 
and there is a
continuous map $f_j:B_j\to B_j$ without periodic points of 
periods less
or equal to $j$. Roughly speaking, $B_j$ was obtained by replacing 
$J_{n\of j}$ in $S_{n\of j}$ 
by a cone over some zero dimensional set $Z_j$.
More precisely, there is a continuous map $q_j$ 
(this map was denoted by
$q$ in \cite{\treelike}) of $B_j$ onto $S_{n\of j}$ 
with the following properties:
\roster
\item"{(Be1)}" $q_j{}^{-1}\of x$ is an one-point set for each 
$x\in \of{S_{n\of j}\setminus
J_{n\of j}}\cup\ofc{d_{n\of j}}$.
\item"{(Be2)}" The set 
$Z_j=q_j{}^{-1}\of{e_{n\of j}}$ 
is zero dimensional.
\item"{(Be3)}"  $q_j{}^{-1}\of{J_{n\of j}}$ 
is a cone over $Z_j$. $q_j{}^{-1}\of{J_{n\of j}}$
is nowhere dense in $B_j$. If $\tilde d$ 
denotes the vertex of the cone and, for each $z\in Z_j$, $A_z$ 
denotes the arc 
between $z$ and $\tilde d$, then $q_j$ restricted to $A_z$ 
is a homeomorphism
onto $J_{n\of j}$.
\item"{(Be4)}" $q_j\circ f_j=g\circ q_j$.
\endroster
It follows that 
\roster
\item"{(Be5)}" every proper subcontinuum of $B_j$ is
arc-wise connected, and
\item"{(Be6)}" if $C$ is a subcontinuum of $B_j$
such that $\tilde d\notin C$, then $C$ is an arc and
$f_j$ restricted to $C$ is a homeomorphism.
\endroster
The continuum $B_j$ cannot be used as $B$ is Section 2.
For instance, not every proper subcontinuum of $B_j$ is an arc. 
To get a continuum that is more suitable for our purpose,  we  
use the technique
presented by J.~B.~Fugate and L.~B.~Mohler in \cite{\fumo}.
As in \cite{\perpoi},
let $\widetilde B_j$ be the inverse limit of the inverse system of 
copies of
$B_j$ with the bonding maps equal to $f_j$. 
Let $\tilde p_j^k$ be the 
projection of $\widetilde B_j$ onto the $k$-th 
element of the inverse system. 
Let $\tilde f_j$ denote the right shift on $\widetilde B_j$, i,e. 
$\tilde f_j\of{\of{b_0,b_1,b_2,\dots}}=
\left(f_j\of{b_0\right),b_0,b_1,b_2,\dots}$.

Let $\pi_j$ denote the map $q_j\circ\tilde p_j^0$. Observe that
the following proposition is true.
\proclaim{Proposition 3.1} $\pi_j{}^{-1}\of x$ is zero-dimensional
for each $x\in S_{n\of j}$.
\endproclaim
\proclaim{Proposition 3.2} Let $A_0$  be an
arc contained in $S_{n\of j}$ 
and let $C$ be a component of $\pi_j{}^{-1}\of{A_0}$.
Then, $\pi_j$ restricted to $C$ is a homeomorphism onto $A_0$.
\endproclaim
\demo{Proof}For each $k=1,2,\dots$, let $A_k=g^{-k}\of{A_0}$.
Since $g$ is a homeomorphism $A_k$ is an arc. By (S2), there
is an integer $m$ such that $d_{n\of j}\notin A_k$ for
each integer $k\ge m$.

By (Be4),
$g^k\circ q_j\circ\tilde p_j^k=\pi_j$ and consequently
$q_j\circ\tilde p_j^k\of C\subset A_k$.
Let $C_k$ be the component of $q_j{}^{-1}\of{A_k}$ containing
$\tilde p_j^k\of C$. Observe that $C$ is the inverse limit of $C_k$'s.

Clearly, $\tilde d\notin C_k$ for each integer
$k\ge m$. 
By (Be6), $C_k$ is an arc and the bonding map
$f_j$ restricted to $C_k$ is a homeomorphism 
for each integer
$k\ge m$.
It follows that $C$ is also
an arc and $\tilde p_j^m$ restricted to $C$ is a homeomorphism 
onto $C_m$. By (Be1) and (Be3), $q_j$ restricted to
$C_m$ is a homeomorphism onto $A_m$. Since $g^m$ restricted
to $A_m$ is a homeomorphism onto $A_0$, we have the result that
$\pi_j=g^m\circ q_j\circ\tilde p_j^m$ restricted to $C$
is a homeomorphism onto $A_0$.
 
\enddemo
\proclaim{Proposition 3.3}Every proper subcontinuum
of $\widetilde B_j$ is an arc.
\endproclaim 
\demo{Proof}
Let $C$ be a proper subcontinuum
of $\widetilde B_j$. Let $C_k=\tilde p_j^k\of C$ for each
integer $k=0,1,\dots$. Since $C_k$ must be a proper
subcontinuum of $B_j$ for
sufficiently large $k$, it follows from (Be5) that
$C_k$ is be a proper
subcontinuum of $B_j$ for each $k$. By (Be5) and
(S1), $A=q_j\of{C_0}$ is an arc. Since 
$C\subset\pi_j{}^{-1}\of{C_0}$, it follows from (3.2) that
$C$ is also an arc.
\enddemo

\proclaim{Proposition 3.4} Let $s$ be a point of
the interval
$\ofb{0,1}$.
Suppose $\kappa$ is a map of
a compactum $M$ onto 
$\ofb{0,1}$ such that
\roster
\item $\kappa^{-1}\of s$ is zero-dimensional ,
and
\item $\kappa$ restricted to each component of $M$ is a 
homeomorphism
onto $\ofb{0,1}$.
\endroster
Then, $M$ is homeomorphic to the product of
$\ofb{0,1}$ and $\kappa^{-1}\of s$. 
\endproclaim
\demo{Proof}
For any $x\in M$, let 
$\sigma\of x$ denote the only point from the component of
$x$ in $M$ such that $\kappa\of{\sigma\of x}=s$. Let
$\varphi:M\to\ofb{0,1}\times\kappa^{-1}\of s$ be defined
by $\varphi\of x=\of{\kappa\of x,\sigma\of x}$. 
We will prove that $\varphi$
is a homeomorphism of $M$ onto $\ofb{0,1}\times\kappa^{-1}\of s$.

Let $x_1,x_2,\dots$ be a sequence of points of $M$ converging
to some point $x\in M$. We will prove that 
$$\lim_{i\to\infty} \sigma\of{x_i}=\sigma\of x.\tag{*}$$
Since $M$ is compact, there is an infinite subsequence 
$\sigma\of{x_{i\of 1}},\sigma\of{x_{i\of 2}},\dots$
 converging to some point
$y\in \kappa^{-1}\of s$.  
Since $M$ is compact, 
$\lim_{k\to\infty}x_{i\of k}=x$ and
$\lim_{k\to\infty}\sigma\of{x_{i\of k}}=y$, $x$ and $y$ must belong
to the same component of $M$. Since $\sigma\of x$ is the only
point of the component of $x$ in $M$ belonging to 
$\kappa^{-1}\of s$, we have the result that $y=\sigma\of x$.
Thus, $\sigma\of x$ is the limit of any convergent subsequence of
$\sigma\of{x_1},\sigma\of{x_2},\dots$, \thetag{*}
is true and consequently
$\sigma:M\to\kappa^{-1}\of s$ is continuous.
It follows that $\varphi$
is also continuous. As an one-to-one continuous 
function defined on a compactum, $\varphi$ is a homeomorphism.
\enddemo

Let $\beta_j:\widetilde B_j\to \ofb{0,1}$ be the map 
$p^0_{n\left(j\right)}\circ q_j\circ\tilde p^0_j$.

\proclaim{Proposition 3.5}
The conditions {\rm (B1)-(B5)} from Section {\rm 2} are satisfied by 
$B=\widetilde B_j$ and $\beta=\beta_j$.
Moreover, if $\lambda$ is the $\beta_j$-length
on $\widetilde B_j$,
then
$\lambda\of{\tilde f_j\of A}=2\lambda\of{A}$ for each arc
$A\subset\widetilde B_j$.
\endproclaim
\demo{Proof}
(B1) follows Proposition 3.3. (B2) follows from Proposition 3.1.
(B3) and (B4) are proven in Proposition 2.11, \cite{\perpoi}.

To prove that (B5) is also satisfied,
suppose $x\in\widetilde B_j$ is not an endpoint. 
By Proposition 3.2, $\pi_j\of x\ne e_{n\of j}$.
By (S4), there is a set $T\subset S_{n\of j}$ containing 
$\pi_j\of x$ in its interior and homeomorphic to the product
of $\ofb{0,1}$ and the Cantor set. By Propositions 3.4 and 3.2,
$\pi_j{}^{-1}\of T$ is homeomorphic to the product
of $\ofb{0,1}$ and the Cantor set. Since $x$ belongs
to the interior of $\pi_j{}^{-1}\of T$, $x$ has a neighborhood
homeomorphic to the product
of the open interval $\of{0,1}$ and the Cantor set.

The remainder of the proposition follows readily from
Proposition 2.12, \cite{\perpoi}.
\enddemo
\demo{Proof of Theorem {\rm 1.1}} 
Let $\lambda$ be the $\beta_j$-length
on $\widetilde B_j$.

By Proposition 2.9, \cite{\perpoi}, the map $\tilde f_j$
does not have
periodic points of periods less than or equal to $j$.
There is a positive number $\eta$ such that any
map $f:\widetilde B_j\to \widetilde B_j$ that is
$\eta$ close to $\tilde f_j$ does have
periodic points of periods less than or equal to $j$.
Chose $\epsilon>0$ so that
the distance between $f$ and $\tilde f_j$ is less than 
$\eta$ for every map 
$f:\widetilde B_j\to \widetilde B_j$ such that
$d_\lambda\of{f,\tilde f_j}<\epsilon$ (see Proposition 2.3).

By Proposition 3.5, we may use Theorem 2.13 with
$B=\widetilde B_j$, $\beta=\beta_j$ and $f_0=\tilde f_j$.
Let $f:\widetilde B_j\to \widetilde B_j$ be the map resulting
from Theorem 2.13. Define $X_j$ as the inverse limit of
copies of $\widetilde B_j$ with $f$ as the bonding map and
let $h_j:X_j\to X_j$ be the right shift map.
Observe that $X_j$ is a tree-like continuum as the inverse limit
of tree-like continua. By Theorem 2.13, 
$X_j$ is hereditarily indecomposable. By the choice of
$\epsilon$, $h_j$ does not have
periodic points of periods less than or equal to $j$.
\enddemo
\Refs
\newcount\refnum
\refnum=0
\def\newref {\global\advance\refnum by 1{\number\refnum}}
{

\ref\no\newref
\by H. Bell
\paper On fixed point properties of plane continua
\jour Trans. Amer. Math. Soc.
\vol 128
\yr 1967
\pages 539--548
\endref

\ref\no\newref
\by D.P.  Bellamy
\paper A tree-like continuum without the fixed-point property
\jour Houston J. Math.
\vol 6 
\yr 1980
\pages1--13
\endref

\ref\no\newref
\bysame
\paper The fixed point property in dimension one 
\inbook in Continua with the Houston problem book
\eds H.~Cook et al.
\publ Marcel Dekker
\publaddr New York
\yr 1995
\pages 27--35
\endref

\ref\no\newref
\by R.~H.~Bing
\paper A homogeneous indecomposable plane continuum
\jour Duke Math. J.
\vol 15 
\yr 1948
\pages 729--742
\endref

\ref\no\newref
\by K. Borsuk
\paper A theorem on fixed points
\jour Bull. Acad. Sci. Polon.\vol 2\yr 1954\pages 17--20\endref

\ref\no\newref
\by L.~Fearnley and D.~G.~Wright
\paper Geometric realization of a Bellamy continuum
\jour Bull. London Math. Soc.
\vol 25 
\yr 1993
\pages 177--183
\endref

\ref\no\newref
\by J. B. Fugate and L. A. Mohler
\paper A note on fixed points in tree--like continua
\jour Topology Proc.
\vol 2 
\yr 1977
\pages 457--460
\endref

\ref\no\newref
\by C. L. Hagopian
\paper A fixed point theorem for plane continua
\jour Bull. Amer. Math. Soc.
\vol  77
\yr 1971
\pages 351--354
\endref

\ref\no\newref
\bysame
\paper Fixed points of plane continua
\jour Rocky Mountain J. Math.
\vol  23
\yr 1993
\pages 119--186
\endref

\ref\no\newref
\bysame
\paper The fixed-point property for simply connected plane continua
\jour Trans. Amer. Math. Soc.
\vol  348
\yr 1996
\pages 4525--4548
\endref

\ref\no\newref
\bysame
\paper The fixed point property for deformations of
tree-like continua
\paperinfo preprint
\endref

\ref\no\newref
\by G. W. Henderson
\paper The pseudo-arc as an inverse limit with one binding map
\jour Duke Math. J.
\vol  31
\yr 1964
\pages 421--425
\endref

\ref\no\newref
\by S. Iliadis
\paper Positions of continua in a plane and fixed points
\jour Vestn. Moskov. Univ.
\vol  4
\yr 1970
\pages 66--70
\endref

\ref\no\newref
\by B. Knaster
\paper Un continu dont tout sous-continu est 
ind\'ecomposable
\jour Fund. Math.
\vol  3
\yr 1922
\pages 247--286
\endref

\ref\no\newref
\by W. Lewis
\paper Continuum theory problems
\jour Topology Proc.
\vol  8
\yr 1983
\pages 361--394
\endref

\ref\no\newref
\by R. Ma\'nka
\paper Association and fixed points
\jour Fund. Math.
\vol  91
\yr 1976
\pages 105--121
\endref

\ref\no\newref
\by R. D. Mauldin (ed.)\book The Scottish Book:
Mathematics from the Scottish Caf\'e
\publ Birkhauser
\publaddr Boston
\yr1981\endref

\ref\no\newref
\by P. Minc
\paper A fixed point theorem for weakly chainable plane continua 
\jour Trans. 	Amer. Math. Soc.
\vol 317 
\yr 1990
\pages 303--312
\endref

\ref\no\newref
\bysame
\paper A tree-like continuum admitting fixed point 
free maps with arbitrarily small trajectories
\jour Topology and its Appl.
\vol 46 
\yr 1992
\pages 99--106
\endref

\ref\no\newref
\bysame
\paper A periodic points free homeomorphism of a tree-like continuum
\jour Trans. Amer. Math. Soc.
\vol 348 
\yr 1996
\pages 1487--1519
\endref

\ref\no\newref
\bysame
\paper A weakly chainable tree-like continuum without the 
fixed point property
\jour Trans. Amer. Math. Soc.
\toappear
\endref

\ref\no\newref
\bysame
\paper A self map of a tree-like continuum with no
invariant indecomposable subcontinuum
\paperinfo preprint
\endref

\ref\no\newref
\bysame and W. R. R. Transue
\paper A transitive map on $\ofb{0,1}$
whose inverse limit is the pseudoarc
\jour Proc. Amer. Math. Soc.
\vol 111 
\yr 1991
\pages 1165--1170
\endref

\ref\no\newref
\by L.G. Oversteegen and J.T. Rogers, Jr.
\paper Tree-like continua as limits of cyclic graphs
\jour Topology  Proc.
\vol  4
\yr 1979
\pages 507--515
\endref

\ref\no\newref
\bysame
\paper An inverse limit description of an atriodic tree-like 
continuum and an induced map without a fixed point
\jour Houston J. Math.
\vol  6
\yr 1980
\pages 549--564
\endref

\ref\no\newref
\bysame
\paper  Fixed-point-free maps on tree-like continua
\jour Topology and its Appl.
\vol  13
\yr 1982
\pages 85--95
\endref

\ref\no\newref
\by K. Sieklucki
\paper On a class of plane acyclic continua with the fixed
point property
\jour Fund. Math.
\vol  63
\yr 1968
\pages 257--278
\endref
}
\endRefs

\enddocument
\end